# Origin of the numerals
# Zero Concept

Ahmed Boucenna
Laboratoire DAC, Department of Physics, Faculty of Sciences,
Ferhat Abbas University 19000 Sétif, Algeria
aboucenna@wissal.dz   aboucenna@yahoo.com

**Abstract**

The partisans of the hypothesis of the Indian origin of the numerals create confusion between the history of the Indian mathematics and the history of our modern numerals. To argue the thesis of the Indian origin of the numbers they confound between: the "intuitive zero" of Brahmagupta, that means "nothing" and which is the difference of two equal numbers, the "numeral zero" used in the representation of the numbers and the "mathematical zero " defined by the modern mathematicians.

"Sifr" designate the ''numeral zero'' and "Shûnya" designate the "intuitive zero". The word "Sifr" is not a traduction of the word "Shûnya" and does not derive from the Indian word "Shûnya", since the word "Sifr" and its derivatives existed in Arabic long before the appearance of zero itself

The facts that the "intuitive zero" and the "mathematical zero" are represented currently by the "numerals zero" symbol "0" are only consequences of the representation of the numbers by the "Ghubari" numerals

**Résumé**.

Les partisans de l'hypothèse de l'origine indienne des chiffres créent une confusion entre l'histoire des mathématiques indiennes et l'histoire de nos chiffres modernes. Pour argumenter la thèse de l'origine indienne des chiffres ils confondent : le ''zéro intuitif''' de Brahmagupta, qui signifie ''rien'' et qui ''est la différence de deux nombres égaux'', le ''chiffre zéro'' utilisée dans la représentation des nombres et le ''zéro mathématique '' définit par les mathématiciens modernes.

"Sifr" désigne le chiffre ''zéro'' et "Shûnya" désigne le ''zéro intuitif''. "Sifr" n 'est pas la traduction de "Shûnya" et ne dérive pas du mot indien "Shûnya" puisque le mot "Sifr" et ses dérivées existaient dans la langue arabes bien avant l'apparition du chiffre ''zéro'' lui même.

Le fait que le ''zéro intuitif '' et le ''zéro mathématique '' soient actuellement représentés par le symbole 0 du ''chiffres zéro'' ne sont que des conséquences de la représentation des nombres par les chiffres "Ghubari".

## 1. Introduction

It is well established that the principle of local value was used by the Babylonians much earlier than by the Hindus and the arabo-islamic [Cantor 1907] and that the Maya of Central America used the principle and symbols for zero in a well-developed numeral system [Bowditch 1910, Morley 1915]. We must recall that our study concerns the origin of the shapes of our ten numerals our numerals and do not concern the Indian mathematics history. The Professor Florian Cajori had well situated the problematic: ''The controversy on the origin of our numerals does not involve the question of the first use of local value and symbols for zero; it concerns itself only with the time and place of the first application of local value to the decimal scale and with the origin of the forms or shapes of our ten numerals.'' [Cajori 1919]. Some assumptions have been put forward as to the origin of the numeral symbols. M. F. Woepcke [Woepcke 1863a, Woepcke 1863b] based essentially on the Medieval Arabic writers from the Eastern Islamic World testimony, considers that the modern numerals as deriving from Indian characters. R. G. Kaye puts in doubt the testimonies of Medieval Arabic Writers from the Eastern Islamic World, look for to invalidate the numeral Hindu origin hypothesis, and try to find an European origin for our numerals [Kaye 1907, Kaye 1909, Kaye 1911a, Kaye 1911b, Kaye 1918, Kaye 1919]. In previous work [Boucenna 2006] we identified the symbols of the "Ghubari" numerals through the mixed pagination of an Arabian Algerian manuscript of the beginning of the 19th century. The correspondence, without ambiguity, between the modern numerals and their eldest the "Ghubari" numerals was be established. The elements of the "Abjadi" calculation ("Hissab El-joummel" or "Guematria"), particularly the notions of the numeric value of the Arabic and Hebrew letters was been recalled. Then relation between the "Ghubari" numeral and the Arabic letter whose numerical value is equal to this numeral was showen [table 1].

Table 1: The "Ghubari" numerals used in the numbering of the sheets
of the manuscript "Kitab khalil bni Ishak El Maliki" [Boucenna 2006]

| Numerical value | Arabic Letter | Likely Initial version of the Ghubari numeral | Ghubari numeral | European version of "Ghubar" Numeral | Final version of the Mashriki numeral | Hebrew letter |
|---|---|---|---|---|---|---|
| 1 | ا | ا | ا | ا | ١ | א |
| 2 | ب | ب | ۲ | 2 | ٢ | ב |
| 3 | ج | ج | ۳ | ۳ | ٣ | ג |
| 4 | ح | ح | ح | 4 | ٤ | ד |
| 5 | ٥ | ٥ | 4 | ح | ٥ | ה |
| 6 | 9 | 9 | 6 | 6 | ٦ | ו |
| 7 | ز | ز | 7 | 7 | ٧ | ז |
| 8 | ح | ح | 8 | 8 | ٨ | ח |
| 9 | ط | ط | 9 | 9 | ٩ | ט |
| 10 | ي | ي |  |  | ٠ | י |
| 90 | ص |  | ص | 0 |  |  |

We also mentioned the permutation of the numbers 4 and 5 in Europe. In old European versions one recovers the order 4 and 5 used in Maghreb and in Spain [Hill 1999], [Kunitzsch 2005].

The partisans of the hypothesis of the Indian origin of our numerals start with demonstrating to you that the Indian knew how to count, that they knew to read and to write the numbers, to the point that one can ask himself if the other civilizations (Egyptians, Greek, Roman, complicates,…) didn't know to do it. They admit very easily, with little required, even contradictory arguments. Thus, on the one hand they affirm that it is Brahmagupta that defines the zero in the VII century and on the other hand they advance that it is in the $V^{th}$ century that the zero donned a specific value: " unknown of the Egyptians, disregarded by the Greeks, denied by the Romans, the zero would have taken its takeoff at your Babylonians, some centuries before our era. Its beginnings were modest: two small hooks representing an empty space. It is in India toward the $V^{th}$ century that it will don a specific value. In 628, the mathematician Brahmagupta defines it as the subtraction of a number by himself. He baptized ''Sunya'' (emptiness) and endowed of a proper existence… " [Institut du monde arabe, 2005].

## 2. The Zero

In fact, it is necessary to distinguish three different types of "zero" : the " intuitive zero ", that means " nothing ", the " numeral zero " used in the representation of the numbers and the " mathematical zero " defines by the modern mathematicians.

### 2.1. The ''intuitive zero''

The " intuitive zero " is roughly equivalent to " nothing ", " rien " and " lashay ". It is known since the oldest times by all societies as primitive either such. If a peasant had 7 bags of wheat initially and that after 10 months he would have consumed 7 bags of it, he knows that he remains to him "zero" bag, that means "nothing". The " intuitive zero " was known by the Egyptians, the Greeks, the Romans, …

This " intuitive zero " meaning " nothing " multiplied or divided by any number gives effectively " nothing ". Even the operation " nothing " divided by " nothing " gives " nothing ", since " nothing " divided by 1000 gives " nothing ", " nothing " divided by 100 gives " nothing "; " nothing " divided by 10 gives " nothing "; " nothing " divided by 1 gives " nothing " and there is not any reason so that " nothing " divided by " nothing " would not give, by a simple intuitive reasoning, " nothing ".

It is this " intuitive zero " that is he " Sunya " of Brahmagupta as his reasoning and his result prove concerning the division of " Sunya " by " Sunya " that gives " Sunya ". It is the division by the " mathematical zero ", unknown by Brahmagupta, witch is empty of sense (impossible). It is not necessary to make to tell Brahmagupta what he didn't say and to assign him the knowledge that he didn't have.

### 2. 2. The numeral zero '' 0 ''

The "numeral zero" is known by the Babylonians (Two small hooks), the Mayan of the Central America [Bowditch 1910, Morley 1915], the Arabo-Islamics (a circle or a point) that baptized it " Sifr " and by our contemporary universal civilization. The " numeral zero " has not been ignored of the Egyptians nor disregarded by the Greeks, nor denied by the Romans. It was to them, simply unknown. Contrary to what the Booklet pretends [Institut du monde arabe, 2005] nothing proves, with certainty that the Indian Brahmagupta. Brahmagupta was knew the " numeral zero " since The earliest reliable record of the use of our numerals with the zero is an inscription of 867 A.D. in India [Cajori 1919].

"Sifr" designate the ''numeral zero'' et "Shûnya" designate the ''intuitive zero''. "Sifr" is not a traduction of "Shûnya". "Sifr" does not derive, from the Indian word "Shûnya" since the word "Sifr" and its derivatives existed in Arabic long before the appearance of zero itself [Ibnou Mandhour LISSAN EL ARAB d-130-136].

The symbol 0 of the numeral zero has been compared to the "sukun" of the Arabic script by the Koranic scholar Abu 'Amr 'Uthman al-Dani, d. 1053 [Kunitzsch 2005, Kunitzsch 2006].

In previous work [Boucenna 2006], we showed that the The Arabo-Islamic civilization represented the number zero by one circle: 0 in Maghreb corresponding to the Arabic letter ص (Sad), first letter of the word Arabian sifr designating the number zero and by one point ּ in the Machrek corresponding to the Hebrew letter י (Yodh), tenth letter of the Hebrew alphabet.

The facts that the " intuitive zero " wich is the difference of two equal numbers, and the " mathematical zero " are represented currently by the symbol 0 of the " numerals zero " are only consequences of the representation of the numbers by the "Ghubari" numerals [Boucenna 2006]. The "mathematical zero " and the "intuitive zero" can pretend to have a philosophical dimension but the choice of the symbols (0) and (ּ) to represent the numeral 0 has not been motivated, as some authors pretend [Ifrah, 1996], [Ifrah, 1999], [Ouaknin, 2004] and [Institut du monde arabe, 2005], by high philosophical considerations.

The discovery of the numeral 0, in its present conception which has simplified the representation of the numbers and the algorithms of the basic operations, was certainly a very crucial event. If the choice of the symbols representing the "Ghubari" numerals 1, 2,…, 9 can be justified by the "Abjadi" numerical values of the Arabic letters chosen to represent them, the symbol representing the "Ghubari" numeral 0 cannot be explained in the same way since the value zero is not an "Abjadi" numerical value. Like all new inventions, it is necessary to choose for it a name and a symbol.

The chosen name is the word "Sifr". Is it an Arabic word? Figure 13 shows a part of an Arabic dictionary page which gives the meaning of the word "Sifr" and its derivatives. One can see that the word "Sifr" is not a foreign word to Arabic. It has not been borrowed from another language to describe a new state, unknown in Arabic, as for example the case of the word "Falsafa" that refers to "philosophy", which is borrowed from Greek to describe a new state. The word "Sifr" does not derive, as some authors,([Ifrah, 1996], [Ifrah, 1999], [Ouaknin, 2004] and [Institut du monde arabe, 2005]) , from the Indian word "Shûnya", since the word "Sifr" and its derivatives existed in Arabic long before the appearance of zero itself [Ibnou Mandhour].

In Arabic the meaning of the word "sifr" describes an empty state which one did not expect. The emptiness described by the word "sifr" is in all cases an abnormal situation. The name "Sifr" given to 0 expresses exactly the role that the 0 must play in the representation of the numbers with the help of the numerals. The numeral 0 must replace a numeral that made defect in a given rank. It does not represent the emptiness, it fills the emptiness.

The chosen symbol to represent the "Sifr" is the symbol 0. Does this symbol derive from an Arabic letter as the other "Ghubari" numerals? If yes, which one? As the value 0 is not an "Abjadi" numerical value, one cannot find an Arabic letter having the good "Abjadi" numerical value, to represent the "Ghubari" numeral 0. The strategy of choice of the symbols must be modified so that one will not take into account the "Abjadi" numerical values of the letters. One could think then about choosing the tenth Arabic letter ي (Yaa) of "Abjadi" numerical value 10 to represent the "Ghubari" numeral 0. This was not the case.

The symbol 0 of the numeral zero has been compared to the "sukun" of the Arabic script by the Koranic scholar Abu 'Amr 'Uthman al-Dani, d. 1053 [Kunitzsch 2005, Kunitzsch 2006]. In my opinion, it is the initial shape of the Arabic letter ص (Sad), first letter of the Arabic word صفر "Sifr", which served to manufacture the symbol of the "Ghubari" numeral 0.

The choice of the symbol 0, has not been motivated therefore as some authors pretend [Ifrah, 1996], [Ifrah, 1999], [Ouaknin, 2004] and [Institut du Monde Arabe 2005] by high philosophical considerations.

## 5. CONCLUSION

The partisans of the hypothesis of the Indian origin of the numerals always create, deliberately confusion and an amalgam between the history of the Indian mathematics and the history of our modern numerals. To argue the thesis of the Indian origin of the numbers they confound voluntarily between: the "zero intuitive" of Brahmagupta : ''Sunya'', that means "nothing", the " numeral zero " used in the representation of the numbers and the "mathematical zero " defined by the modern mathematicians.

"Sifr" designate the ''numeral zero'' et "Shûnya" designate the ''intuitive zero'', thus "Sifr" is not a traduction of "Shûnya". "Sifr" does not derive, from the Indian word "Shûnya" since the word "Sifr" and its derivatives existed in Arabic long before the appearance of zero itself

The facts that the " intuitive zero " wich is the difference of two equal numbers, and the " mathematical zero " are represented currently by the symbol 0 of the " numerals zero " are only consequences of the representation of the numbers by the "Ghubari" numerals [Boucenna 2006].